\documentclass[11pt]{amsart}
\usepackage[english]{babel}

\usepackage{amsmath, amssymb}
\usepackage[mathscr]{eucal}
\usepackage{amsfonts,amscd}
\usepackage{amsthm,latexsym}
\usepackage{fullpage}
\usepackage{bm}
\usepackage{graphicx}
\usepackage{enumitem}
\newtheorem{thm}{Theorem}

\def\Z{\mathbb Z}

\def\R{\mathbb R}

\def\F{\mathbb F}
\def\O{\operatorname{O}}

\def\Tr{\operatorname{\text{ Tr }}}

\def\mod{\rm $mod \,$ }
\date\today
\title{The error term in the Sato-Tate theorem of Birch}

\author[M. Ram Murty]{M. Ram Murty}
\address{M. Ram Murty, Queen's University, Kingston, Ontario- K7L 3N6,  Canada  }
\email{murty@queensu.ca}
\author[Neha Prabhu]{Neha Prabhu}
\address{Neha Prabhu, Queen's University, Kingston, Ontario- K7L 3N6,  Canada  }
\email{neha.prabhu@queensu.ca}
\thanks{}

\begin{document}
	
	\begin{abstract}
		We establish an error term in the Sato-Tate theorem of Birch. That is, for $p$ prime, $q=p^r$   we show that $$\#\{ (a,b) \in \F_q^2 : \theta_{a,b}\in I\} =\mu_{ST}(I)q^2 + \O_r(q^{7/4})$$ for any interval $I\subseteq[0,\pi]$ where for an elliptic curve $E: y^2= x^3 +ax +b$, the quantity $\theta_{a,b}$ is defined by $2\sqrt{q}\cos\theta_{a,b} = q+1-E(\F_q) $ and $\mu_{ST}(I)$ denotes the Sato-Tate measure of the interval $I$. 
	\end{abstract}
	
	
	\maketitle
	
	\section{Introduction} In 1968, Birch \cite{Birch} proved that the Sato-Tate conjecture holds for the family of elliptic curves
	 \begin{equation*}
	 	y^2= x^3 + ax + b \bmod p
	\end{equation*}
as $(a,b)$ varies over elements in $(\Z / p\Z)^2$ such that $4a^3 + 27b^2 \not\equiv 0\,  \mod p$, where $p$ is a fixed prime. Let $(\frac{\cdot}{p})$ denote the Legendre symbol. More precisely, he proved that \begin{equation*}
\sum\limits_{a,b=0}^{p-1} \left| \sum_{x=0}^{p-1} \left( \frac{x^3+ax+b}{p}\right)\right|^{2R} \sim \frac{1}{R+1}\binom{2R}{R}p^{R+2}
\end{equation*}
as $p \to\infty$. These are the moments predicted by the Sato-Tate conjecture and by standard probability theory one can deduce the relevant distribution. This moment calculation then implies the Sato-Tate distribution for the angles $\theta_{a,b}$ where we write 
\begin{equation*}
a_p(a,b) = -\sum_{x=0}^{p-1}\left( \frac{x^3+ax+b}{p}\right) = 2\sqrt{p}\cos\theta_{a,b}
\end{equation*}
where, if $E(\F_p)$ denotes the number of $\F_p-$ points (including the point at infinity) on an elliptic curve $y^2= x^3 +ax +b$, then  $a_p(a,b)$ is $p+1- E(\F_p)$.
There were two key ingredients in Birch's proof. The first was Deuring's theorem \cite{Deuring} that there are a total of $H(t^2-4p)$ isomorphism classes of elliptic curves over $\F_p$ with $p+ 1 -t$ points where $H(N)$ denotes the Hurwitz-Kronecker class number. The second ingredient is the Eichler-Selberg trace formula (see \cite[Appendix]{Lang}) which gives the trace of the $n$-th Hecke operator acting on the space of holomorphic cusp forms of even weight $k\geq 4$ with respect to the full modular group. 

A natural question that arises from Birch's paper is the order of the error term in the Sato-Tate distribution. More precisely, let us fix an interval $I= [\alpha, \beta] \subseteq [0,\pi]$. We want to count the number $$N_I(p): =\# \left\{ (a,b) \in (\Z/p \Z)^2 : p\nmid \Delta(a,b), {\theta_{a,b}} \in I \right\}$$ and give a quantitative estimate for \begin{equation} \label{discrepancy}
| N_I(p) -\mu_{ST}(I)p^2| \end{equation} where the Sato-Tate measure of the interval $I$ is given by $\mu_{ST}(I)= \frac{2}{\pi}\int_I \sin^2( \theta)d\theta. $

Using the discrepancy estimate of Niederreiter \cite{Neid}, Banks and Shparlinski \cite{BSh} noted that \begin{equation}\label{Birch error bound}
N_I(p) -\mu_{ST}(I)p^2 = \O(p^{7/4})
\end{equation}
would follow from the work of Katz \cite{Katz} extended to their setting. This extension is not routine and appears in the work of Michel \cite{Mic} where he deals with the case of one-parameter families of elliptic curves. There is also a related paper of Fisher \cite{Fisher}.
The work of Michel \cite{Mic} relies heavily on Weil II of Deligne \cite{Del2}. For many of us working in classical analytic number theory, Weil II  and its cohomological mysteries present formidable prerequisites that often represent a ``black box" whose pronouncements must be accepted on faith. On the other hand, using the moment estimates in Birch's proof, \eqref{discrepancy} was estimated by Miller and Murty in \cite{Miller-Murty} as well, where only a log saving over the trivial estimate of $p^2$ was obtained. 
The goal of this paper is to show that the estimate \eqref{Birch error bound} can be deduced using classical techniques, from just the Ramanujan-Petersson conjecture, now a theorem due to Deligne (which is implied by Weil I \cite{Del1} and \cite{Del3}) as well as the two key ingredients of Birch \cite{Birch} mentioned earlier. The result is in fact, true in the more general case of elliptic curves over a finite field $\F_q=\F_{p^r}$. Let $$N_I(q) := \# \left\{ (a,b) \in \F_q^2 : \Delta(a,b) \neq 0,\, {\theta_{a,b}} \in I \right\},$$ where we now have $a_q(a,b) = q+1 -E(\F_q) = 2\sqrt{q}\cos\theta_{a,b}.$ We prove the following.
\begin{thm}\label{main}
	Assume the notation above. We then have
	\begin{equation}
	N_I(q) -\mu_{ST}(I)q^2 = \O_r(q^{7/4}).
	\end{equation}
\end{thm}
We note that the result of Banks and Shparlinski, (and that in Theorem \ref{main}) gives a true error term only when the size of the interval $I$ is greater that $p^{-1/4 +\varepsilon}$. This was improved (on average) by Baier and Zhao in \cite{BZ} and recently by David, Koukoulopoulos and Smith \cite{DKS} where the effective version of Birch's theorem is shown to hold for intervals $I$ of length $\geq p^{-1/2 +\varepsilon}$ although in these cases, the saving is only a power of a logarithm over the main term. 

\subsection*{Acknowledgements} The authors would like to thank Igor Shparlinski for his helpful comments and the anonymous referee for suggestions that improved the exposition of this article.

\section{Preliminaries}
\subsection{Isomorphism classes of elliptic curves}\label{Isomorphism classes}
We briefly discuss the ingredients from the theory of counting elliptic curves, which will be needed for the proof. For more details, see \cite{Lenstra} or \cite{Schoof}. 

For $p \neq 2,3$, consider the elliptic curve over $\F_q$ in Weierstrass form
$$E_{a,b}: y^2= x^3 +ax +b.$$
Analogous to the case of $\F_p$, the number of $\F_q$ points on $E$, given by $E(\F_q)$ is $q + 1 - a_q(a,b)$. We have Hasse's bound $|a_q(a,b)|\leq 2\sqrt{q}$. 
Two curves $E=E_{a,b}$ and $E'=E_{a',b'}$ over $\F_q$ are isomorphic if there is an element $u\in\F_q^{*}$ such that $a'=u^4a$ and $b'=u^6b$. An automorphism  of $E$ is an isomorphism from $E$ to $E$. Clearly, isomorphism of elliptic curves is an equivalence relation and the size of the equivalence class of $E$ is given by $\frac{\# \F_q^{*}}{\#\text{Aut } E}$. 

For $t^2 < 4q$, Deuring \cite{Deuring} essentially showed that the number isomorphism classes of elliptic curves $E$ with $ q+ 1 -t$ points, weighted by $\#(\text{Aut }E)^{-1}$ is $H(t^2-4q)$, where $H(N)$ is the Hurwitz-Kronecker class number (see \cite{Lenstra} for a detailed description of these numbers). Thus, for $t^2<4q$, the total number of curves $E$ over $\F_q$ with $q+1-t$ points is $(q-1)H(t^2-4q)$. 
\subsection{Chebyshev polynomials}\label{Chebyshev}
The Chebyshev polynomials of the second kind, $U_n(x)$ for integers $n\geq 0$ are defined recursively in the following way:
\begin{equation*}
\begin{split}
U_0(x) &=1 \\
U_1(x) &= 2x\\
U_n(x) &= 2xU_{n-1}(x) - U_{n-2}(x), \quad n\geq 2.
\end{split}
\end{equation*} 
If $x= \cos \theta$, then the polynomials can be written explicitly as $$U_n(x) = \frac{\sin(n+1)\theta}{\sin \theta}.$$
In our application of these polynomials,
$x= \cos\theta_{a,b}$. It is not hard to see that 
 \begin{equation*}
 \frac{\rho^{k-1} -\bar{\rho}^{k-1}}{\rho-\bar{\rho}}= q^{\frac{k-2}{2}}U_{k-2}\left(\frac{t}{2\sqrt{q}}\right) 
\end{equation*}
where $\rho$ and $\bar{\rho}$ are the roots of the equation $y^2 -ty + q =0$. 
Observe that $$2\cos(n\theta) = \frac{\sin(n+1)\theta}{\sin\theta} - \frac{\sin(n-1)\theta}{\sin\theta} = U_{n}(\cos \theta) -U_{n-2}(\cos\theta),$$ a fact that will be needed later.
\subsection{Beurling-Selberg Polynomials}
The Beurling-Selberg polynomials have been frequently used to obtain effective results on equidistribution and in this note, we use this method to study the quantity $N_I(p)$. We give a brief introduction to these polynomials, see \cite[Chapter 1]{Mont} for a detailed exposition.
Let $\chi_J(x)$ denote the characteristic function of the interval $J = [\alpha,\beta] \subseteq [0,1]$ and $M \geq 1$ be an integer.  One can construct trigonometric polynomials $S_{M,J}^-(x)$ and $S_{M,J}^+(x)$ of degree less than or equal to $M,$ respectively called the minorant and majorant Beurling-Selberg polynomials for the interval $J$ such that
\begin{itemize}
	\item[(a)] For all $x \in \R,\,S_{M,J}^-(x)\leq \chi_I(x) \leq S_{M,J}^+(x)$
	\item[(b)] $$\widehat{S}^{\pm}_{M,J}(0)= \int_{0}^{1} S_{M,J}^{\pm}(x) dx = \beta - \alpha \pm \frac{1}{M+1},$$
	\item[(c)] For $ 0 < |m| \leq M,$
	\begin{equation*}\label{FC}
	\left|\widehat{S}_{M,J}^{\pm}(m) - \widehat{\chi}_J(m)\right| \leq \frac{1}{M+1}.
	\end{equation*}
\end{itemize}
Henceforth, we will supress the $J$ in the subscript of the Fourier coefficients with the understanding that the definition of these approximating polynomials depends on the interval $J$. 

 For $J= [\alpha, \beta]$ we will also use the following estimates, which follow from properties (b) and (c) listed above.
\begin{equation}\label{S(0)}
\widehat{S}^{\pm}_M(0)= (\beta -\alpha) \pm \frac{1}{M+1}
\end{equation} and for $m>0,$ 
\begin{equation}\label{S(m)}
\widehat{S}^{\pm}_{M}(m) + \widehat{S}^{\pm}_{M}(-m) = \frac{\sin2\pi m \beta - \sin 2\pi m \alpha}{m\pi} + \O\left(\frac{1}{M+1}\right).
\end{equation}
These polynomials were used in \cite{MS} to study the `vertical' Sato-Tate distribution in the case of modular forms.
\vspace{0.5cm}

\section{Proof of the main theorem}
Let $I \subseteq [0,\pi]$. We consider the angles $ \left\{ \frac{\theta_{a,b}}{2\pi}, - \frac{\theta_{a,b}}{2\pi}  \right\}$ and count when they occur in $I'= I/2\pi$. 
Approximating using Beurling-Selberg polynomials, we have
\begin{equation*}
\begin{split}
N_I(q) =\sum_{{a,b\in \F_q}\atop { \Delta(a,b)\neq 0}} \chi_{I'}\left(\frac{\theta_{a,b}}{2\pi}\right) &\leq \sum_{{a,b\in \F_q}\atop { \Delta(a,b)\neq 0}} \left(\sum_{|m|\leq M}\widehat{S}^+_{M}(m)e\left(m\frac{\theta_{a,b}}{2\pi}\right) + \sum_{|m|\leq M}\widehat{S}^+_{M}(m)e\left(-m\frac{\theta_{a,b}}{2\pi}\right)\right)\\
&= \sum_{{a,b\in \F_q}\atop { \Delta(a,b)\neq 0}}\sum_{|m|\leq M}\widehat{S}^+_{M}(m)2\cos\left(m\theta_{a,b}\right)\\
&= 2\widehat{S}^+_{M}(0)\sum_{{a,b\in \F_q}\atop { \Delta(a,b)\neq 0}}1 + \sum_{0<|m|\leq M}  \widehat{S}^+_{M}(m)\sum_{{a,b\in \F_q}\atop { \Delta(a,b)\neq 0}} \left[ U_{m}(\cos\theta_{a,b}) - U_{m-2}(\cos\theta_{a,b}) \right],
\end{split}
\end{equation*}
where, as noted in Section \ref{Chebyshev}, $$U_{r}(\cos\theta_{a,b}) = \dfrac{\sin(r+1)\theta_{a,b}}{\sin\theta_{a,b}}$$ denotes the $r$-th Chebyshev polynomial of the second kind evaluated at $\cos(\theta_{a,b})$. 
 Note that $U_0(x)=1$.  Therefore, 
\begin{equation}\label{approximation}
\begin{split}
N_I(q) &\leq q^2 \left( 2\widehat{S}^+_{M}(0)-(\widehat{S}^{\pm}_{M}(2) + \widehat{S}^{\pm}_{M}(-2)) \right) + \sum_{m=1}^{2} (\widehat{S}^{\pm}_{M}(m) + \widehat{S}^{\pm}_{M}(-m))\sum_{{a,b\in \F_q}\atop { \Delta(a,b)\neq 0}} U_m(\cos\theta_{a,b}) \\ & \qquad + \sum_{3\leq m\leq M} (\widehat{S}^{\pm}_{M}(m) + \widehat{S}^{\pm}_{M}(-m))\sum_{{a,b\in \F_q}\atop { \Delta(a,b)\neq 0}} \left[ U_{m}(\cos\theta_{a,b}) - U_{m-2}(\cos\theta_{a,b}) \right].
\end{split}
\end{equation}

Using \eqref{S(0)} and \eqref{S(m)}, we see that 
\begin{equation}\label{main-term}
2\widehat{S}^{\pm}_{M}(0)-(\widehat{S}^{\pm}_{M}(2) + \widehat{S}^{\pm}_{M}(-2)) = \mu_{ST}(I) + \O\left(\frac{1}{M+1}\right)
\end{equation} where  $\mu_{ST}(I) = \frac{2}{\pi}\int_I \sin^2 \theta d\theta $ denotes the Sato-Tate measure of $I\subseteq [0,\pi]$. It remains to estimate the sums $$\sum_{{a,b\in \F_q}\atop { \Delta(a,b)\neq 0}} U_{m}(\cos\theta_{a,b}) = \sum_{{a,b\in \F_q}\atop { \Delta(a,b)\neq 0}} \dfrac{\sin(m+1)\theta_{a,b}}{\sin\theta_{a,b}}$$ for $m=1,\ldots,M$.  
If we write $a_q(a,b)= t=2\sqrt{q}\cos\theta_{t}$ where $|t|\leq 2\sqrt{q}$, then we have \begin{equation}\label{Chebyshev (a,b) to (t,p)}
\sum_{{a,b\in \F_q}\atop { \Delta(a,b)\neq 0}} \dfrac{\sin(m+1)\theta_{a,b}}{\sin\theta_{a,b}} =(q-1) \sum_{|t|\leq2\sqrt{q}}H(t^2-4q)\frac{\sin(m+1)\theta_{t}}{\sin\theta_{t}},
\end{equation}
where we group the curves into isomorphism classes, as discussed in Section \ref{Isomorphism classes}. If $m$ is odd, writing $$\frac{\sin(m+1)\theta_{t}}{\sin\theta_{t}} = U_m\left(\frac{t}{2\sqrt{q}}\right),$$ we see that the sum in \eqref{Chebyshev (a,b) to (t,p)} is zero. This follows from the fact that in the polynomial $U_m\left(\frac{t}{2\sqrt{q}}\right),$ the parity of the powers of $t$ that appear is the same as that of $m$, so the terms corresponding to $t$ and $-t$ cancel each other when $m$ is odd.

On the other hand, using the Eichler-Selberg Trace formula (see \cite{Lang} or \cite{KnLi}) we have for even $k\geq 4$,
\begin{equation*}
\begin{split}
\Tr T_k(q) &= \frac{k-1}{12}q^{\frac{k-2}{2}}\delta(q,2) - \frac{1}{2} \sum\limits_{|t|\leq 2\sqrt{q}} \frac{\rho^{k-1} -\bar{\rho}^{k-1}}{\rho-\bar{\rho}} H(t^2-4q) - \frac{1}{2}\sum_{dd'=q}\min(d,d')^{k-1}\\
&= \frac{k-1}{12}q^{\frac{k-2}{2}}\delta(q,2) -\frac{1}{2}\sum_{|t|\leq2\sqrt{q}}q^{\frac{k-2}{2}} \frac{\sin(k-1)\theta_{t}}{\sin\theta_{t}} H(t^2-4q) - \frac{1}{2}\sum_{dd'=q}\min(d,d')^{k-1},
\end{split}
\end{equation*} where $\delta(q,2)$ is $1$ when $q$ is a square and zero otherwise. 
 Using the Ramanujan-Petersson bound for Hecke eigenvalues, we have 
 \begin{equation} \label{trace}
\Tr T_k(q) \ll kq^{\frac{k-1}{2}}, 
 \end{equation} 
 using the fact that the dimension of the space of cusp forms of weight $k$ and full level grows like $k$. Therefore, 
\begin{equation} \label{Hurwitz class number}
\sum_{|t|\leq2\sqrt{q}}H(t^2-4q)\frac{\sin(k-1)\theta_{t}}{\sin\theta_{t}} \ll rkq^{1/2}, 
\end{equation}
 where $q=p^{r}$. 
 Going back to \eqref{Chebyshev (a,b) to (t,p)}, letting $k= m+2$, we deduce that for even $m=2,\ldots,M$,
\begin{equation}\label{Chebyshev-sum-bound}
\sum_{{a,b\in \F_q}\atop { \Delta(a,b)\neq 0}} \dfrac{\sin(m+1)\theta_{a,b}}{\sin\theta_{a,b}} \ll rmq^{\frac{3}{2}}.
\end{equation}
Using \eqref{main-term}, \eqref{Chebyshev-sum-bound} and the estimate $(\widehat{S}^{\pm}_{M}(m) + \widehat{S}^{\pm}_{M}(-m)) \ll 1/m$ (evident from \eqref{S(m)}) in \eqref{approximation}, we get 
\begin{equation}\label{final error}
N_I(q)-q^2\mu_{ST}(I) \ll \frac{q^2}{M} + rMq^{3/2}.
\end{equation}
Letting $M= \left\lfloor p^{1/4}r^{-1/2} \right\rfloor$, we get $$N_I(q)-q^2\mu_{ST}(I) \ll_r q^{7/4}.$$

Using $\widehat{S}^-_{M}(m)$ the lower bound estimation is similar.

\section{Concluding remarks}
It is interesting to consider to what extent our error term is best possible. For example, for $q=p$, the sum in \eqref{trace} is essentially $p^{-\frac{k-1}{2}} \Tr T_k(p)$. If we accept the prediction that the Sato-Tate distribution corresponding to distinct Hecke eigenforms as discussed in \cite{Murty}, then we can arrange all the Fourier coefficients appearing in Tr $T_k(p)$ to be arbitrarily close to $2p^{\frac{k-1}{2}}$ simultaneously, for infinitely many primes $p$. Thus the estimate in \eqref{Hurwitz class number} cannot be improved for all primes $p$. This does not however say anything about the combined error term in \eqref{final error}. Thus, the question of the optimal error term becomes an intruiging problem for future research.

\end{document}